\documentclass{amsart}
\usepackage{hyperref}
\usepackage{graphicx}
\usepackage{amscd}
\usepackage{amsmath}
\usepackage{amssymb}
\usepackage{verbatim}
\newtheorem{conjecture}{Conjecture}

\usepackage[dvipsnames]{xcolor}
\usepackage{pgfplots}
\pgfplotsset{compat=newest}
\usetikzlibrary{arrows.meta}
\usetikzlibrary{bending}
\usepackage{bm}
\usepackage{asymptote}
\newtheorem{theorem}{Theorem}[section]

\theoremstyle{definition}
\newtheorem{definition}[theorem]{Definition}
\newtheorem{proposition}[theorem]{Proposition}
\newtheorem{example}[theorem]{Example}

\theoremstyle{remark}
\newtheorem{remark}[theorem]{Remark}

\numberwithin{equation}{section}

\begin{document}

\title{Lyapunov Exponents of variations of Hodge structures with $G_2$ monodromy}

\author{Genival Da Silva Jr.}

\address{\newline Department of Mathematics \newline Eastern Illinois University\newline Charleston, IL}
\email{jrrieman@gmail.com}


\begin{abstract}
We investigate the Lyapunov Exponents of a variation of Hodge structure which has $G_2$ as geometric monodromy group, and discuss formulas for the sum of positive Lyapunov Exponents of variations of Hodge structures of any weight.
\end{abstract}

\maketitle

\section{Introduction}
The main object of study in this paper is the dynamics of a weight $2$ variation of Hodge structures $\mathcal{V}=(\mathbb{H},V,\nabla)$ coming from the transcendental cohomology of a family of Elliptic surfaces over $\mathbb{P}^1-\{0,\pm \frac{2}{3\sqrt{3}},\infty\}$, with hodge numbers $(2,3,2)$. The monodromy group is a subgroup of $G_2(\mathbb{Z})$. An explicit description of the monodromy matrices and further discussion can be found in \cite{jr1}.

In order to study the Lyapunov Exponents of the variation above, we follow the ideas from \cite{filip1}, that is, we use the representation theory of $G_2$ to study the exponents. We first recall some background and notation.

\begin{theorem}(Oseledets)
Let $(X,\mu)$ denotes a probability measure space, equipped with an ergodic flow $g_t$. Suppose $V\to X$ is a vector bundle, equipped with a norm $||\cdot ||$ and a lift of $g_t$. The lift being an action on the fibers of $V$:
$$g_t(x):V_x\to V_{g_t x}$$
If the linear maps above are integrable:
\begin{equation}\label{integrable}
\int_X \sup_{t\in[-1,1]} ||g_t(x)||d\mu(x)<\infty
\end{equation}
Then there are numbers $\lambda_1 > \lambda_2 > \ldots > \lambda_n$, called \textit{Lyapunov Exponents}, and a measureable $g_t$-invariant decomposition
$$V_x=\bigoplus_i V_x^{\lambda_i}$$
such that for $\mu$-a.e. $x\in X$ and $v\in V_x^{\lambda_k}$ we have:
$$\lim_{t\to\pm\infty} \frac{1}{t}\log||g_t(x)v||=\lambda_k $$
\end{theorem}
For the sake of convenience, we write the Lyapunov Exponents as a list $\lambda_1 \geq \lambda_2 \geq \ldots \geq \lambda_d$, where $d=\dim V_x$. That is, we repeat the exponent $\lambda_i$ a number of times equal to the dimension of the Oseledets space $V_x^{\lambda_i}$. We also refer to $\lambda_1 \geq \lambda_2 \geq \ldots \geq \lambda_d$ as the \textit{Lyapunov spectrum} of the cocycle $g_t$.
\begin{definition}
A Hodge structure of weight $w$ consists of an abelian group $H_\mathbb{Z}$ together with a decomposition of its complexification by vector subspaces:
$$H_\mathbb{Z}\otimes \mathbb{C} = \bigoplus_{p+q=w} H^{p,q}$$
such that $\overline{H^{p,q}}=H^{q,p}$. Given $H^{p,q}$ we can define a decreasing filtration $F^p=\bigoplus_{k\geq p}H^{k,w-k}$, and conversely, given $F^p$, we can define $H^{p,q}=F^p\cap \overline{F^q}$. We say that the Hodge structure is \textit{polarized} if it comes equipped with a non-degenerate $(-1)^w$-symmetric bilinear form $Q(x,y)$ such that $Q(F^p,F^{w-p+1})=0$ and the form $i^{p-q}Q(x,\bar{y})$ is positive definite.
\end{definition}
\begin{definition}
A \textit{variation of Hodge structures} of weight $w$ over a complex manifold $B$ consists of the following data:
\begin{itemize}
    \item[$\circ$] A local system $\mathbb{H}$ over $B$
    \item[$\circ$] A decreasing filtration of the vector bundle $V:= \mathbb{H}_\mathbb{C}\otimes \mathcal{O}_B$ by vector sub-bundles
    $$V=F^0\supseteq \ldots \supseteq F^{n-1} \supseteq F^n$$
    such that $F^p_x$ induces a Hodge structure of weight $w$ on the fiber $(\mathbb{H}_\mathbb{C})_x$ for each $x\in B$
    \item[$\circ$] The natural flat (Gauss-Manin) connection $\nabla$ on $V$ satisfies
    $$\nabla F^p \subseteq F^{p-1}\otimes \Omega^1_B$$
\end{itemize}
The variation is said to be polarized if in addition we have a bilinear map $Q(x,y)$ on $\mathbb{H}$ which induces a polarization on each fiber of $\mathbb{H}_\mathbb{C}$.
\end{definition}
Many examples of variations of Hodge structures come from geometry, that is, from the cohomology of a family of projective algebraic varieties. 
\begin{example}(Legendre Family)
Let $B:=\mathbb{P}^1-\{0,1,\infty\}$ and for $t\in B$, let $X_t\subset \mathbb{P}^2$ be the Elliptic curve:
$$y=x(x-1)(x-t)$$
Set $\mathcal{X}=\{ (x,t)\in \mathbb{P}^2\times B | x\in X_t \}$ and consider the natural morphism:
$$\pi:\mathcal{X}\to B$$
whose fiber at $t\in B$ is $X_t$. Then the local system $\mathbb{H}=R^1\pi_*\mathbb{Z}$ underlies a variation of Hodge structures of weight $1$.
\end{example}
\begin{figure}
    \begin{tikzpicture}[scale = .9]
      \begin{axis}[
        axis lines = middle,
        domain = 0 :4,
        xlabel = $x$,
        ylabel = $y$,
        xtick = {0},
        ytick = {0},
        samples = 100,
        thick,
        declare function = {
          H0(\x) = (\x * (\x-1)*(\x-1))^0.5;
          H1(\x) = -(\x * (\x-1)*(\x-1))^0.5;
          H2(\x) = (\x * (\x-1)*(\x-2))^0.5;
          H3(\x) = -(\x * (\x-1)*(\x-2))^0.5;
          H4(\x) = (\x * (\x-1)*(\x-3))^0.5;
          H5(\x,\y) = (\x * (\x-1)*(\x-\y))^0.5;
        }, ]
        \addplot[domain=0:3,blue] {H0(x)}; 
        \addplot[domain=0:3,blue] {H1(x)};
        \addplot[domain=0:2,cyan] {H2(x)}; 
        \addplot[domain=0:2,cyan] {H3(x)}; 
        \addplot[domain=2:3,cyan] {H2(x)}; 
        \addplot[domain=2:3,cyan] {H3(x)}; 
        \addplot[domain=0:2,cyan] {H4(x)}; 
        \addplot[domain=0:2,cyan] {-H4(x)}; 
        \addplot[domain=1.7:3,cyan] {H5(x,1.7)}; 
        \addplot[domain=1.7:3,cyan] {-H5(x,1.7)}; 
        \addplot[domain=1.5:3,cyan] {H5(x,1.5)}; 
        \addplot[domain=1.5:3,cyan] {-H5(x,1.5)}; 
        \addplot[domain=1.3:3,cyan] {H5(x,1.3)}; 
        \addplot[domain=1.3:3,cyan] {-H5(x,1.3)}; 
        \addplot[domain=2.3:3,cyan] {H5(x,2.3)}; 
        \addplot[domain=2.3:3,cyan] {-H5(x,2.3)}; 
        \addplot[domain=2.6:3,cyan] {H5(x,2.6)}; 
        \addplot[domain=2.6:3,cyan] {-H5(x,2.6)}; 
        \addplot[domain=0:0.7,cyan] {H5(x,0.7)}; 
        \addplot[domain=0:0.7,cyan] {-H5(x,0.7)}; 
        \addplot[domain=1:3,cyan] {H5(x,0.2)}; 
        \addplot[domain=1:3,cyan] {-H5(x,0.2)}; 
      \end{axis}
    \end{tikzpicture}
    \label{legendre}
    \caption{The legendre family plot in $\mathbb{R}^2$}
\end{figure}
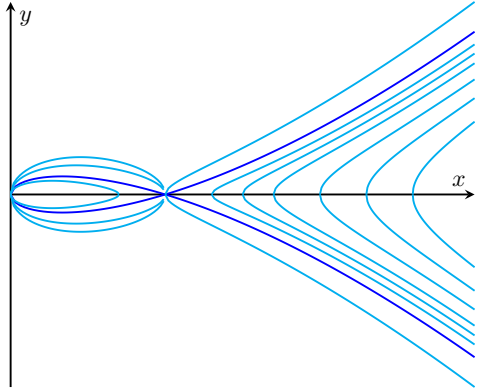
Using the notation above, recall that to give a local system over $B$ underlying a polarized variation of Hodge structure is the same as giving a representation $\rho:\pi_1(B,x_0)\to Aut(V_{x_0},Q)$. The image of this representation is called the monodromy group.
\begin{example}
The example above has $Aut(V_{x_0},Q)=Sp(2,\mathbb{Z})=SL(2,\mathbb{Z})$ and the monodromy group is the congruence group $\Gamma(2)$ of index 6.
\end{example}
\begin{remark}
Except for a few cases, there's no algorithm to explicitly compute the monodormy group of any variation. Even the computation of the index inside $Aut(V_{x_0},Q)$ can be quite challenging. Interestingly, the Kontsevich formula tends to be valid in cases where the monodromy group has infinite index.
\end{remark} 

Suppose that $B$ is a $1$ dimensional complex manifold having a metric of constant negative curvature $-4$, thus we view $B$ as a hyperbolic Riemann Surface. Let $g_t$ be the unitary geodesic flow over $B$. Using the Gauss-Manin connection $\nabla$, we can lift the flow $g_t$ to a flow $G_t$ acting on $V_x$, for each $x\in B$. Hence, $G_t$ defines a cocycle and we are interested at its Lyapunov exponents.
\begin{remark}
The fact that $G_t$ satisfies equation \ref{integrable} is not obvious but follows from the curvature properties of period domains. \cite[Proposition 3.2]{eskin1}
\end{remark}
We want to explore the fact that the monodromy group in our case is a Semisimple Lie Group, so we need a version of the Oseledets theorem that takes advantage of this fact. The approach here is the same as in \cite{filip1}.

Consider a principal bundle $P\to X$ with a right $G$-action. Given a real representation $\rho: G\to GL_n(\mathbb{R})$, we define the associated bundle to $\rho$ as $W=P\times_G\mathbb{R}^n\to X$, where $P\times_G\mathbb{R}^n=\{ (p,v)\in P\times\mathbb{R}^n | (p,v)\cong (pg,\rho(g^{-1})v)\;\forall g\in G\}$.

Let $g_t$ be an ergodic flow acting on a probability space $(X,\mu)$. Suppose the flow lifts to a principal $G$-bundle $P\to X$, with $G$ a semisimple Lie Group, such that the lift satisfy the integrability condition \ref{integrable}.

Let $\Phi$ be the root system of $G_\mathbb{C}$, $\Phi_{res}$ the root system of the real form $G$, and $\mathfrak{a}$ be the lie algebra of a maximal split torus of $G$. Denote by $\mathfrak{a}_+$ the positive Weyl chamber of $\mathfrak{a}$ and $r:\Phi\to\Phi_{res}$ the restriction map.

\begin{theorem}(Kaimanovich)\cite{kai,filip1}\label{kaima}
There exists a Lyapunov vector $\Gamma\in\mathfrak{a}_+$ which controls Lyapunov exponents as follows. Given a representation $\rho:G\to Aut(V)$ with weights $\Sigma$ and associated bundle $W=P\times_G V$. The Lyapunov exponents of the flow $g_t$ on $W$ are given by evaluating the restricted weights $r(\Sigma)$ on $\Gamma$.
\end{theorem}

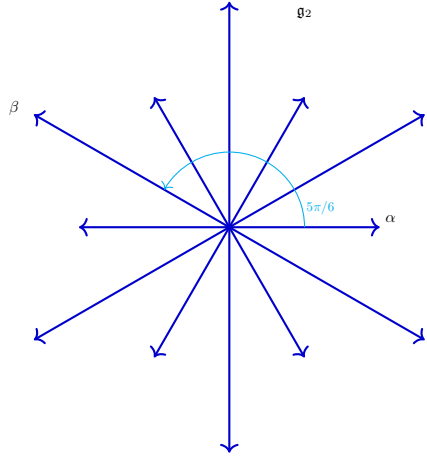
\begin{figure}
    \centering
    \begin{tikzpicture}
    \foreach\ang in {60,120,...,360}{
     \draw[->,blue!80!black,thick] (0,0) -- (\ang:2cm);
    }
    \foreach\ang in {30,90,...,330}{
     \draw[->,blue!80!black,thick] (0,0) -- (\ang:3cm);
    }
    \draw[cyan,->](1,0) arc(0:150:1cm)node[pos=0.1,right,scale=0.5]{$5\pi/6$};
    \node[anchor=south west,scale=0.6] at (2,0) {$\alpha$};
    \node[anchor=south west,scale=0.6] at (-3,1.414) {$\beta$};
    \node[anchor=north,scale=0.6] at (1,3) {$\mathfrak{g}_2$};
 \end{tikzpicture}
    \caption{The root system of $\mathfrak{g}_2$}
    \label{root_system}
\end{figure}

\section{Fundamental representations of $\mathfrak{g}_2$}

Recall that the root system of the lie algebra $\mathfrak{g}_2$ has type $G_2$ and is composed of $12$ roots, with $2$ being simple, as depicted in figure $\ref{root_system}$. 
Let $e_1,e_2,e_3$ be the standard basis of $\mathbb{R}^3$. Then
the roots are
$$\pm (e_1-e_2),\pm (e_1-e_3), \pm (e_2-e_3), \pm (-e_1+2e_2-e_3), \pm (2e_1-e_2-e_3), \pm (e_1+e_2-2e_3)$$
and one choice of simple roots is:
$$\alpha_1=e_1-e_2,\; \alpha_2=-e_1+2e_2-e_3$$
The restricted root system of ${\mathfrak{g}_2}(\mathbb{R})$ is also of type $G_2$, with simple roots:
$$\beta_1=f_1-f_2,\; \beta_2=-f_1+2f_2-f_3$$
The restriction map $r$ is given by (see \cite{vinberg})
$$r(\alpha_1)=\beta_1,\;r(\alpha_2)=\beta_2$$
Therefore, we also have
\begin{equation*}
    \begin{split}
        r(e_i-e_j) &= f_i-f_j\\
        r(2e_1-e_2-e_3) &= 2f_1-f_2-f_3\\
        r(e_1+e_2-2e_3) &= f_1+f_2-2f_3
    \end{split}
\end{equation*}
The lie algebra $\mathfrak{g}_2$ has $2$ fundamental representations. The first one $\rho_{std}$, called \textit{standard}, is a $7$ dimensional representation with weights
$$\{0,\pm (e_1-e_2),\pm (e_1-e_3), \pm (e_2-e_3)\}$$
and restricted weights
$$\{0,\pm (f_1-f_2),\pm (f_1-f_3), \pm (f_2-f_3)\}$$
The second one $\rho_{adj}$ is the $14$ dimensional \textit{adjoint} representation, the weights are the roots plus $0$ with multiplicity two, the dimension of the cartan subalgebra.

Now using theorem \ref{kaima} and the discussion above we have
\begin{proposition}
Let $g_t$ be an ergodic flow acting on a probability space $(X,\mu)$. Suppose that the flow $g_t$ lifts to a principal $G_2$-bundle $P\to X$ such that the lift satisfies the integrability condition \ref{integrable}. Let $W_{std}$ be the associated vector bundle to the representation $\rho_{std}$, with structure group $G_2$. Then the Lyapunov Exponents of $(g_t,W_{std})$ are:
$$0,\pm (\lambda_1-\lambda_2),\pm (\lambda_1-\lambda_3), \pm (\lambda_2-\lambda_3)$$
\end{proposition}
\section{Example}
An example of a variation which has monodromy group $\mathbb{Q}$-Zariski dense in $G_2$ can be found as follows.

Let $\mathcal{E}\rightarrow\mathbb{P}^1:y^2=x(x-1)(x-z^2)$ be a rational elliptic surface with singular fibers at $z=-1,0,1,\infty$. For $t\neq 0,\pm \frac{2}{3\sqrt{3}} , \infty$, take a base change by:
\begin{equation}
    E_t\rightarrow \mathbb{P}^1:w^2=tz(z-1)(z+1)+t^2
\end{equation}
The result is a family of elliptic surfaces $X_t\rightarrow E_t$ with 7 singular fibers on each surface, as described below:
\begin{equation}
    \begin{array}{ccc}
        \mathcal{X} & \hookleftarrow & X_t\\
        \downarrow \pi &  & \downarrow \pi_t\\
        \mathcal{E} & \hookleftarrow & E_t\\
        \downarrow &  & \downarrow\\
        \mathbb{P}^1 & \hookleftarrow & \{ t \}
    \end{array}
\end{equation}
Then as $t\neq 0,\pm \frac{2}{3\sqrt{3}}$ varies, $H^2_{tr}(X_t)$ define a variation of Hodge structure $\mathcal{H}^2$ of weight 2 and geometric monodromy group $G_2$, see \cite{jr1}. The monodromies are:
\begin{equation}
    \begin{split}
        M_{\frac{2}{3\sqrt{3}}} &= \begin{bmatrix}
            1&0&1&0&0&0&0\\
            0&1&0&1&0&0&0\\
            0&0&1&0&0&0&0\\
            0&0&0&1&0&0&0\\
            0&0&0&0&1&0&0\\
            0&0&0&0&0&1&0\\
            0&0&0&0&0&0&1
            \end{bmatrix}\\
            \\
        M_{-\frac{2}{3\sqrt{3}}} &= \begin{bmatrix}
            1&0&0&0&0&0&0\\
            0&1&0&0&0&0&0\\
            -1&0&1&0&0&0&0\\
            0&-1&0&1&0&0&0\\
            0&0&0&0&1&0&0\\
            0&0&0&0&0&1&0\\
            0&0&0&0&0&0&1
            \end{bmatrix}\\
            \\
        M_0 &= \begin{bmatrix}
            1&2&-2&-2&2&-1&-2\\
            -2&-3&6&2&-4&3&6\\
            2&6&-3&-2&6&-3&-4\\
            -2&-2&2&1&-2&1&2\\
            0&0&-4&0&1&-2&-4\\
            -4&-4&4&4&-4&1&4\\
            0&-4&0&0&-4&2&1
            \end{bmatrix}\\
            \\
        M_\infty &= \begin{bmatrix}
            0&-4&1&0&-4&2&2\\
            4&0&4&1&-2&2&4\\
            -1&4&-3&-2&6&-3&-4\\
            0&-1&2&1&-2&1&2\\
            -4&0&-4&0&1&-2&-4\\
            0&0&4&4&-4&1&4\\
            0&-4&0&0&-4&2&1
            \end{bmatrix}
    \end{split}
\end{equation}
In this example, rank of $\mathcal{H}^2$ is 7 and the intersection form has signature $(4,3)$. Therefore, the Lyapunov Spectrum has shape:
$$\lambda_1\geq\lambda_2\geq\lambda_3\geq 0\geq -\lambda_3\geq-\lambda_2\geq-\lambda_1$$
We also could have predicted the shape of the spectrum by using the proposition above. Additionally, by \cite[Theorem 1]{eskin-matheus1}, there's actually a spectral gap, i.e. distinct exponents:
$$\lambda_1>\lambda_2>\lambda_3$$
Using the same ideas and notation of \cite[section 4.2]{filip1}, a formula for $\lambda_1$ can be found, but there's a technicality involved in this example. The Hodge bundle $\mathcal{H}^{2,0}$ has dimension $2$ in our case, in particular is not a line bundle, so the formula
$$\lambda_1=\int \left(\lim_{t\to \infty}\frac{1}{2\pi t}\int_0^{2\pi} \log |g_t^\theta v|d\theta\right)d\eta(v)$$
is not necessarily independent of the $v$ taken, and will not simplify as much as we had hoped.

\section{Sum of positive Lyapunov exponents}
For variations of Hodge structures of weight $1$, we have the celebrated \textit{Kontsevich formula} for the sum of positive Lyapunov Exponents.
\begin{theorem}(Kontsevich)\cite{eskin2}
Let $\mathcal{H}^1$ be a variation of Hodge structures of weight $1$ and rank $2g$ over an hyperbolic Riemann surface $C$ with singular set $S$ of the local system. Then the sum of positive Lyapunov exponents associated to this variation satisfy:
$$\lambda_1+\ldots+\lambda_g=\frac{2\int_C c_1(\mathcal{H}^{1,0})}{2g(C)-2+\#S}$$
\end{theorem}

For weights higher than $1$ the formula is not always true, that is, for weight $n$ the sum is not necessarily $\int c_1(F^{\lceil \frac{n}{2}\rceil})$, where $\lceil . \rceil$ is the ceil function. This was first shown by Kontsevich in the case of weights 3 variations coming from the cohomology of Calabi-Yau threefolds. On the other hand, He also found cases where the formula works in the same setting of weight 3 as well.

So a natural question is when the formula is valid and if not valid, what adjustments can be made in order for it to work. Kontsevich himself suggested that the answer may have to do with the arithmeticity of the monodromy group inside its $\mathbb{Q}$-Zariski closure. More precisely, he found that the cases where his formula worked had one thing in common, the monodromy group was thin, i.e not a finite index subgroup, whereas in the cases where it failed, the monodromy group was an arithmetic subgroup. It is not known whether or not thin monodromy implies Kontsevich formula in general.

The question of arithmeticity of monodromy groups is highly non trivial and there is no universal approach to prove whether or not a given monodromy group is arithmetic or not, see \cite{sarnak}. Apart from a few cases \cite{brav}, the problem remains unsolved.

Motivated by the result of Kontsevich, we have the following:
\begin{conjecture}
Let $\mathcal{H}^n$ be a weight $n$ variation of Hodge structures with thin monodromy, over a hyperbolic Riemann surface $C$, and singular locus $S$. Suppose $\lambda_1,\ldots,\lambda_k$ are the positive Lyapunov exponents. If $k=\dim F^{\lceil \frac{n}{2}\rceil}$, then their sum satisfies:
$$\lambda_1+\ldots+\lambda_k=\frac{2\int_C c_1(F^{\lceil \frac{n}{2}\rceil})}{2g(C)-2+\#S}$$
If $k<\dim F^{\lceil \frac{n}{2}\rceil}$ but $k>\dim \mathcal{H}^{n,0}$, set $d:=\dim \mathcal{H}^{n,0}$, then:
$$\lambda_1+\ldots+\lambda_d=\frac{2\int_C c_1(\mathcal{H}^{n,0})}{2g(C)-2+\#S}$$
\end{conjecture}
In particular, if $d=1$ the conjecture simplifies to:
\begin{proposition}
If the conjecture above is true and $d=1$, then the top Lyapunov exponent is rational and is equal to
$$\lambda_1 = \frac{2 deg(\mathcal{H}^{n,0})}{2g(C)-2+\#S}$$
\end{proposition}
A general Kontsevich formula, whether the monodromy is thin or not, can be obtained by expanding the expression
\begin{equation}\label{variational}
    \partial\bar\partial \log |\bigwedge c_i|^2
\end{equation}
where $c_i$ are a basis of an isotropic subspace of $(\mathcal{H}^n_\mathbb{R})_x$ for some $x\in C$, and $|.|$ is the Hodge norm, see \cite[Appendix B]{filip2} for the weight 1 case. Using the definition of Lyapunov Exponents and Green's theorem for Laplacian, the computation of the sum of positive exponents boils down to simplify equation \ref{variational} using well known formulas of complex geometry.

The difficulty in obtaining a general formula and simple expression is due to the fact that the Hodge bundles $F^{\lceil \frac{n}{2}\rceil}$ or $\mathcal{H}^{n,0}$ are rarely a line bundle in general, so the formulas for the curvature and connection don't simplify as they do in weight 1 case. Nevertheless, we believe a closed expression can be found and we expect to say more about this in future works.

\subsection*{Acknowledgements} I thank Carlos Matheus for very helpful discussions during the preparation of this paper.

\bibliographystyle{amsplain}

\end{document}